\documentclass[11pt]{amsart}
\usepackage{amscd}
\usepackage{amssymb}
\usepackage{boxedminipage}
\usepackage{enumerate}
\usepackage[all]{xy}
\usepackage[utf8]{inputenc}
\usepackage[T1]{fontenc}

\newtheorem{thm}{Theorem}[section]

\newtheorem{cor}[thm]{Corollary}
\newtheorem{conjecture}[thm]{Conjecture}

\newtheorem{question}[thm]{Question}
\theoremstyle{definition}
\newtheorem{defn}[thm]{Definition}

\newtheorem{remark}[thm]{Remark}

\def\i{^{-1}}

\def\mod{\operatorname{mod}}
\def\dim{\operatorname{dim}}
\def\ker{\operatorname{ker}}

\def\Hilb{\operatorname{Hilb}}

\def\|{\mid}

\def\iso{\cong}

\def\Ext{\operatorname{Ext}}

\def\Hom{\operatorname{Hom}}

\def\Hilb{\operatorname{Hilb}}

\begin{document}
\title [Quotients of Koszul algebras with almost linear resolution]
{Qoutients of Koszul algebras with almost linear resolution}
\author {Jon Eivind Vatne}
\address{Department of mathematics, physics and computer science\\
Faculty of Engineering\\
Bergen University College\\
P.O.Box 7030\\
N-5020 Bergen\\
Norway}   
        
\email{jev@hib.no }

\begin{abstract}
We compute the Hilbert series, and the graded vector space structure, of Ext-algebras of quotients of Koszul algebras with almost linear resolution.
The example of the generic determinantal varieties is treated in detail.
\end{abstract}

\maketitle

\section{Introduction}
The Koszul condition for associative algebras was introduced by Priddy \cite{Prid}.
This condition is satisfied by many interesting algebras, and one can often perform explicit computations on algebras with this property.
The condition has been generalized in various ways, including to operads \cite{GK}, to algebras with relations in a fixed higher degree \cite{B01}, and to algebras with relations in any degree \cite{CS}.
Examples of computations in these cases can be found in many articles, including \cite{BF}, \cite{AS87}, \cite{PP}, \cite{FV06}, \cite{BG}, \cite{GMMVZ}, \cite{Val}, \cite{Bar}, \cite{GKR}.
See also \cite{Lod1},\cite{Lod2}.\\

In this article, we will see that it is possible to compute the graded vector space structure on Ext-groups of quotients of Koszul algebras, whose resolution is linear except for the first map, which is homogeneous of a fixed degree $d\geq 4$.
We say that the resolution is {\em almost linear}.
The Ext-groups are isomorphic to graded pieces of the tensor product of the Koszul dual and a tensor algebra based on the resolution, see Theorem \ref{main}.

If $A$ is A Koszul algebra, $R$ the quotient with almost linear resolution, we prove this result by ``reverse engineering'' on the spectral sequence
\[E^{pq}_2=\Ext_R^p(k,\Ext_A^q(R,k))\Rightarrow \Ext_A^{p+q}(k,k).\]
There is one drawback; we have to assume that the degree $d\neq 3$ for our proof to work.  Since the case $d=2$ is already known, we will throughout assume that $d\geq 4$.\\

Examples can be found for instance among determinantal varieties, and these examples point towards a possible notion of $N-$Koszul commutative algebra.

\section{Preliminaries on Koszulity}
We collect a few definitions and results that will be needed later on.

\subsection{Koszul algebras and Koszul duality}
A graded algebra
\[A=\bigoplus_{m\geq 0} A_m\]
with $A_0=k$, the ground field, and generated by the finite dimensional vector space $V=A_1$  is called {\em Koszul} if the minimal resolution of the ground field is linear:
\[ \cdots \rightarrow A(-2)^{a_2}\rightarrow A(-1)^{a_1}\rightarrow A\rightarrow k\]

Under these conditions, we can write $A$ as a quotient of a tensor algebra
\[A=T(A_1)/I,\]
where the ideal $I$ is generated in degree 2, so that
\[I_2\subset A_1\otimes A_1.\]
The Koszul dual algebra is
\[A^!=T(A_1^\ast)/(I_2^\perp),\]
where $V^\ast=\Hom_k(V,k)$ is the vector space dual, and
\[I_2^\perp =\ker(A_1^\ast\otimes A_1^\ast \rightarrow I_2^\ast).\]
When $A$ is Koszul, so is $A^!$.  Furthermore,
\[A_m^!\iso \Ext^m_A(k,k),\]
and the algebra structure on $A$ is the same as the Yoneda algebra structure on Ext.
From the grading of $A$, it is clear that the Ext-groups have an internal degree, and in the Koszul case the internal degree of $A_m^!$ is concentrated in degree $-m$ (this is a way to rephrase the linearity).

\subsection{Generalizations of the Koszul condition}

R. Berger introduced the notion of N-Koszul algebras in \cite{B01}.  A graded algebra as above, except that the ideal $I$ is generated in degree $N$, is called $N$-Koszul if the graded resolution of the ground field is alternately linear and of degree $N-1$:
\[\cdots \rightarrow A(-2N)^{a_4}\rightarrow A(-N-1)^{a_3} \rightarrow A(-N)^{a_2}\rightarrow A(-1)^{a_1}\rightarrow A\rightarrow k\]

Then the internal grading of the Ext-algebra $\Ext^*_A(k,k)$ places $\Ext_A^i(k,k)$ in a single degree, $-2i$ if $i$ is even, $-2i-1$ if $i$ is odd.  In this case, the $A_\infty$-structure on the Ext-algebra, has that the compositions $m_2$ and $m_N$ are the only nonzero compositions, se \cite{FV06}.\\

T. Cassidy and B. Shelton \cite{CS} introduced a more general condition for Koszulity, namely that the associative structure on the Ext-algebra should be generated by $\Ext^1$ and $\Ext^2$.
With the restrictions on degrees of ideal generators for Koszul and $N-$Koszul algebras, this definition exactly captures the Koszul and the $N-$Koszul condition, see for instance \cite{GMMVZ}.\\

The Koszul condition, under its various guises, makes certain computations feasible.  See for example \cite{AS87}, \cite{FV06}, \cite{PP}, \cite{BG}, \cite{BF}.

\section{Quotients of Koszul algebras}

We fix a Koszul algebra $A$, with a quotient algebra $R$.

\begin{defn}
The resolution of $R$ as a (left) $A$-module is
\[\cdots \rightarrow P_n\rightarrow P_{n-1}\rightarrow \cdots \rightarrow P_1\rightarrow P_0\]
where $P_0=A$.
If the ideal $I=\ker(A\rightarrow R)$ has all its minimal generators in the same degree $d$, and that apart from the first map $P_1\rightarrow P_0$, all differentials are {\em linear}, we say that $R$ has an {\em almost linear resolution.}
Note that this forces the resolution to be minimal.
Thus
\[P_i\iso A(-d-i+1)^{b_i} \text{ for } i\geq 1,\, P_0\iso A.\]
\end{defn}

In this situation, we pick out the generating graded vector space $G_i$ for the $i$th syzygy, so that
\[P_i=G_i \otimes_k A.\] 
Thus
\[G_i\iso k(-d-i+1)^{b_i}\]
in its internal grading, but it also has cohomological grading placing it in degree $i+1$.  The tensor algebra $T(\oplus G_i^\vee)$, where $()^\vee$ indicates vector space dual, is graded using the cohomological grading with $G_i^\vee$ in degree $i+1$, but also inherits an internal grading from the internal grading of the $G_i^\vee$.
This grading convention will make sense because of the following theorem, and its proof:

\begin{thm}\label{main}
Assume that $R$ is a quotient of a Koszul algebra $A$ with an almost linear resolution, and relations in degree $d\geq 4$.
Then
\[\Ext^m_R(k,k)\iso (A^!\otimes T(\oplus G_i^\vee))_m\]
as (internally) graded vector spaces, 
where $A^!$ is the Koszul dual of $A$ and the $G_i$ are the generating graded vector spaces of the syzygies.
\end{thm}

\begin{proof}
We consider the change of rings spectral sequence for $\Ext$ groups:
\[E^{pq}_2=\Ext_R^p(k,\Ext_A^q(R,k))\Rightarrow \Ext_A^{p+q}(k,k).\]
We know the inner Ext-groups from the projective resolution of $R$;
\[\Ext_A^\ast(R,k)\iso\Hom_A(P_\ast,k)\iso G_\ast^\vee.\]
Thus
\begin{equation}\label{produktstruktur}
\Ext_R^\ast(k,\Ext_A^\ast(R,k))\iso\Ext_R^\ast(k, G_\ast^\vee)\iso\Ext_R^\ast(k,k)\otimes G_\ast^\vee.
\end{equation}
In words, the $i$th column in the $E^2$-term is the zeroeth column times the zeroeth row term in the $i$th column.\\

Steps in the proof:
First we consider what is known about the limit and the inital data of this spectral sequence.  Our theorem will be proved by identifying the zeroeth row ($q=0$) of the $E_2$-term, since
\[E_2^{p0}=\Ext_R^p(k,\Ext_A^0(R,k))\iso\Ext_R^p(k,k).\]
Second: We show that the differentials $d_m$ are zero whenever the image is not on the zeroeth row, and are injective if the image is on the zeroeth row.  We will prove this by using the internal degree only.\\

Third: The results about the differentials imply a recurrence relation on the $\Ext^p_R(k,k)$ which can be solved to give the theorem.\\

First we look at the limit term.
By Koszulity the spectral sequence converges to $A^!$, and the graded pieces of the limit term are concentrated on the zeroeth row of the spectral sequence.
More precisely, $E_\infty^{p0}=A^!_p=\Ext_A^p(k,k)$, with internal degree $-p$.
Since our spectral sequence is in the first quadrant, this limit is achieved for $E_{p+1}$, though the spectral sequence as a whole does not degenerate at any finite step.\\

Then we look at the initial conditions.
The assumption on the resolution of $R$ means that the zeroeth column ($p=0$) in the $E^2$-term is known,
\[E_2^{0q}\iso G_q^\vee\iso k(-q-d+1)^{b_q},\,q>0.\]

Since all generators for the ideal of $R$ in $A$ have degree $d\geq 4$, 
\[\Ext_A^1(k,k)\iso\Ext_R^1(k,k)\iso k(-1)^{b_1}.\]
Since there aren't any differentials going into the first column, we get one copy of the zeroeth column, twisted by $(-1)$, for each summand $k(-1)$ in the first Ext-group.
In other words,
\[E_2^{1q}\iso \bigoplus k(-q-d),\,q>0.\]

At this point, we now the two leftmost parts of the $E_2$-term.  Writing only the internal degrees, we have that
\[\begin{array}{cccc}
\vdots & \vdots & \vdots & \\
-d-2 & -d-3 & ? &\cdots\\
-d-1 & -d-2 & ? &\cdots\\
-d & -d-1 & ? &\cdots\\
0 & -1 & ? & \cdots
\end{array}\]

We now come to the second stage of the proof.
Consider first the map

\[d_2:E_2^{01}\rightarrow E_2^{20}.\]
This is the only differential going into position $(2,\,0)$, and so the cokernel must be the limit term $A^!_2$.  Also, it is the only differential going out of position $(0,\,1)$, so it is necessarily injective:
\[0\rightarrow E_2^{01} \rightarrow E_2^{20} \rightarrow A^!_2 \rightarrow 0.\]
The internal degree of the leftmost term is $-d$, of the rightmost $-2$.  Thus the internal degrees of $E_2^{20}$ are $-d$ and $-2$.\\

{\em We will show that the internal degrees at any given point in the spectral sequence are constant modulo d-2}.\\

For the two leftmost columns, this is true (there is only one degree at each point), and at position $(2,\,0)$ the internal degrees are $-d$ and $-2$, both congruent to $-2$ modulo $d-2$.  Using Equation \ref{produktstruktur}, we see that
\[E_2^{2q} \text{ has degrees } -d-1-q \text{ and } -2d-q+1,\,\,q\geq 1,\]
so the degrees are constant modulo $d-2$ also in this case, congruent to $-q-3$.

Writing only the internal degrees modulo $d-2$, we now have this image of the three leftmost columns of the $E_2$-term:
\[\begin{array}{ccccc}
\vdots & \vdots & \vdots & \vdots & \\
-d-2 & -d-3 & -d-4,-2d-2 & ? &\cdots\\
-d-1 & -d-2 & -d-3,-2d-1& ? &\cdots\\
-d & -d-1 & -d-2,-2d&  ? &\cdots\\
0 & -1 & -2,-d & ? & \cdots
\end{array}\sim \begin{array}{ccccc}
\vdots & \vdots & \vdots & \vdots & \\
-4 & -5 & -6 & ? &\cdots\\
-3 & -4 & -5 & ? &\cdots\\
-2 & -3 & -4 &  ? &\cdots\\
0 & -1 &  -2 & ? & \cdots
\end{array}\]
At this point, it is clear that the differential $d_2:E_2^{0,q}\rightarrow E_2^{2,\,q-1}$ must be zero for $q>1$, since the source has degree $-q-1$ and the image has degree $-2-q$ modulo $d-2$.\\

We now proceed by induction on the column number $p$, showing three statements at once:
\begin{itemize}
\item[1)] The internal degrees of $E_2^{p0}$ are all congruent to $-p$ modulo $d-2$, the internal degrees of $E_2^{pq}$ are all congruent to $-p-q-1$ modulo $d-2$ for $q\geq 1$ .
\item[2)] The differentials $d_s:E_s^{p-s,s-1}\rightarrow E_s^{p0}$ are injective.
\item[3)] The differentials $d_s:E_s^{p-s,s-1+q}\rightarrow E_s^{pq}$ are zero for $q\geq 1$.
\end{itemize}
The start of the induction is already taken care of; we have seen these properties for $p\leq 2$.\\

Assume that all these properties hold to the left of the $p$-th column.
The possible degrees for $E_2^{p0}$ are $-p$, coming from the $A^!_p$ in the limit, and the degree coming from differentials going into $E_2^{p0}$, namely the degrees of
\[E_2^{p-2,\,1},\,E_2^{p-3,\,2},\,\cdots,\,E_2^{p-k,k-1},\,E_2^{p-p,p-1}=E_2^{0,p-1}.\]
By induction, each of these has internal degrees congruent to
\[-(p-k)-(k-1)-1\cong -p \mod (d-2).\]
Therefore all internal degrees of $E_2^{p0}$ are congruent to $-p$ modulo $d-2$, and by Equation \ref{produktstruktur}, $E_2^{pq}$ has internal degrees congruent to
\[-q-d+1-p\cong-p-q-1 \mod (d-2).\]
1) is therefore true also for column $p$.\\

For property 3), note that the internal degrees of $E_s^{pq}$ must be a subset of the internal degrees of $E_2^{pq}$, which are all congruent to $-p-q-1$ modulo $d-2$ when $q\geq 1$.
The source of the differential $d_s:E_s^{p-s,s-1+q}\rightarrow E_s^{pq}$ likewise has internal degrees congruent to
\[-(p-s)-(s-1+q)-1= -p-q,\]
different from $-p-q-1$ modulo $d-2$ (recall that $d\geq 4$).
Therefore each differential $d_s$ is zero when the target $E_2^{pq}$ has $q\geq 1$, and 3) is proved.\\

Finally, 2) follows since the differential $d_s:E_s^{p-s,s-1}\rightarrow E_s^{p0}$ is the only one starting from position $(p-s,\,s-1)$ that is not zero (by induction), and the limit term $E_\infty^{p-s,\,s-1}$ is zero, so the differential must be injective.  This concludes the inductive proof of the three statements.\\

In order to identify to prove that the terms 
\[E_2^{p0}\iso \Ext_R^p(k,k)\iso (A^!\otimes T(\oplus G_i^\vee))_p\]
we now use property 2) inductively, to get an equality of graded vector spaces:
\[\begin{array}{rl}
E_2^{p0} \iso & A^!_p\oplus E_2^{p-2,1} \oplus E_2^{p-3,2}\oplus \cdots \oplus E_2^{0,p-1}\\
\iso & A^!_p\oplus (E_2^{p-2,0}\otimes G_1^\vee)\oplus (E_2^{p-3,0}\otimes G_2^\vee)\oplus \cdots \oplus (E_2^{00}\otimes G_{p-1}^\vee)\\
\iso & A^!_p\oplus (A^!\otimes T(\oplus G_i^\vee))_{p-2}\otimes G_1^\vee \oplus (A^!\otimes T(\oplus G_i^\vee))_{p-3}\otimes G_2^\vee \oplus \\
& \oplus \cdots \oplus (A^!\otimes T(\oplus G_i^\vee))_0\otimes G_{p-1}^\vee
\end{array}\]
Since $G_i^\vee$ is in cohomological degree $i+1$, this proves the theorem.
\end{proof}

\begin{cor}\label{maincor}
Under the assumptions of the theorem, the Hilbert series of $\Ext_R^\ast(k,k)$ is the product of the Hilbert series of $A^!$ and $T(\oplus G_i^\vee)$:
\[\Hilb(\Ext_R^\ast(k,k))(t)=\Hilb(A^!)(t)\Hilb(T(\oplus G_i^\vee))(t)=\frac{\Hilb(A^!)(t)}{1-b_1t^2-b_2t^3-b_3t^4-\cdots}.\]
\end{cor}
\begin{remark}
The internal degree can be captured by introducing a second variable $u$, so that we get a two-variable Hilbert series
\[\Hilb(\Ext_R^\ast(k,k))(t,u)=\sum \dim_k\Ext_R^i(k,k)_{-s}t^iu^s,\]
and then the corollary can be refined to
\[\Hilb(\Ext_R^\ast(k,k))(t,u)=\frac{\Hilb(A^!)(tu)}{1-b_1u^{d}t^2-b_2u^{d+1}t^3-b_3u^{d+2}t^4-\cdots}.\]
\end{remark}

\begin{remark}[d=3]
Part 1) of the inductive argument is still true for $d=3$, but being congruent modulo $3-2=1$ is not interesting.
In fact, there will be equal degrees in the source and target of differentials many places in the spectral sequence in case $d=3$, so parts 2) and 3) of the inductive argument can not be proved by looking only on the degrees.\\

On the other hand, the statement of the theorem was found partly by computer experiments in Macaulay2, \cite{M2}, mostly in the case $d=3$.
\end{remark}

\begin{conjecture}
The theorem holds also in case $d=3$.
\end{conjecture}

\begin{remark}[d=2]
The theorem holds for $d=2$, and also the proof, if equality modulo $d-2$ is exchanged by {\em equality}.
This reproves a classical result, see for instance \cite{PP} Chapter 2.5.
\end{remark}

\section{Example: generic determinantal varieties}
The symmetric algebra $S(V)$ is Koszul, so we can in particular discuss qoutients of this algebra with almost linear resolutions, and their associated projective schemes.\\

Let $M$ be a matrix of homogeneous polynomials,
\[M=(f_{ij}), \, i=1,\cdots ,n,\,j=1,\cdots ,m.\]
Assume $4\leq n\leq m$.
Let $I$ be the ideal of maximal ({\em i.e.} $n\times n$) minors of $M$.
Consider in particular the case where 
\[M=(x_{ij})\]
is a matrix of variables, so that $I$ is the ideal of the generic deteminantal variety in $\mathbb{P}^{nm-1}=\mathbb{P}(V)$, and let $R=k[x_{ij}]/I$ be the homogeneous coordinate ring.\\

The purpose of this section is to present the following example computation using the main theorem.
\begin{thm}
The Ext-algebra of the homogeneous coordinate ring $R$ of the generic determinantal variety in $\mathbb{P}^{nm-1}=\mathbb{P}(V)$ ($4\leq n\leq m$) satisfies
\[\Ext^l_R(k,k)\iso (\Lambda V^\ast \otimes T(\oplus G_i^\vee))_l\]
where $\Lambda V^\ast$ is the exterior algebra and $G_{i+1}$ is the tensor product of the (dual of the) $i$th symmetric power of an $n$-dimensional space and the $(i+n)$th exterior power of an $m$-dimensional space.
The Hilbert series is
\[\Hilb \Ext^\ast_R(k,k)(t) = \frac{(1+t)^{nm}}{1-\binom{m}{n}t^2-n\binom{m}{n+1}t^3-\cdots - \binom{m-2}{m-n-1}\binom{m}{m-1}t^{m-1}-\binom{m-1}{m-n}t^m}.\]
\end{thm}
\begin{proof}

To prove this, we consider the matrix $M$ as a map between free $S(V)$-modules of ranks $m$ and $n$.
Then our algebra $R$ can be considered the cokernel of the $n$th exterior power of this map, and this is resolved by the {\em Eagon-Northcott complex}, see for example Eisenbud's \cite{Eis} Theorem A2.10.
In the same reference, an explicit description of the terms of the Eagon-Northcott complex is given, with $i$th term equal to the tensor product of the (dual of the) $i$th symmetric power of a rank $n$ module and the $(i+n)$th exterior power of a rank $m$ module.
Eisenbud resolves the ideal, not the quotient algebra, so we must change the index by one, and take the generating graded vector space as before.
The differentials in the Eagon-Northcott complex are all linear, so the conditions of Theorem \ref{main} are fulfilled, and the first equality in the theorem follows.\\

The formula for the Hilbert series similarly follows from Corollary \ref{maincor}, by recalling that the Hilbert series of the exterior algebra is
\[\Hilb \Lambda V^\ast=(1+t)^{\dim_k V}=(1+t)^{mn}.\]
\end{proof}

The same proof would give a similar result for other matrices of linear forms, as long as the Eagon-Northcott complex is a resolution.
In particular, consider 
\[M=\left(\begin{array}{ccccccc}
x_0 & x_1 & \cdots & x_s & 0 & \cdots & 0\\
0 & x_0 & \cdots & x_{s-1} & x_s & \cdots & 0\\
\vdots & \vdots & & \vdots & \vdots & & \vdots\\
0 & \cdots & x_0 & \cdots & \cdots & x_{s-1} & x_s\\
\end{array}\right),\]
in which case $I$ id the $n$-th power of the irrelevant ideal in $k[x_0,\,\cdots ,x_s]$.
Here $m=n+s$.
In this case, we get that the Hilbert series of the Ext-algebra of the Artin algebra $R=k[x_0,\,\cdots ,x_s]/(x_0,\,x_2,\cdots ,x_s)^n$ is
\[\Hilb \Ext_R^{ast}(k,k)(t) =\frac{(1+t)^s}{1-\binom{n+s}{n}t^2-n\binom{n+s}{n+1}t^3-\cdots - \binom{n+s-2}{s-1}\binom{n+s}{n+s-1}t^{n+s-1}-\binom{n+s-1}{s}t^{n+s}}.\]

See \cite{GKR} for a similar result for the rational normal curve.

\section{Context for the result}

The generalization of the Koszul condition due to Cassidy and Shelton is often easy to check negatively, by considering degrees of Ext-groups (see \cite{CS} page 100, where it is shown that the four-dimensional Artin-Schelter algebras with two generators will not satisfy this condition for degree reasons).
The quotient algebras of Koszul algebras with almost linear resolution pass this test, {\em i.e.} the degrees satisfy the same properties as the degrees of generalized Koszul algebras in this sense.

\begin{question}
Are quotients of Koszul algebras with almost linear resolutions generalized Koszul algebras in the sense of Cassidy and Shelton?
More generally, what is the algebra structure of the Ext-algebra?
\end{question}

The example of the generic determinantal variety studies a commutative algebra, all of whose relations are in a fixed degree.
There is an operadic notion of Koszul algebra over a Koszul operad, and a commutative algebra is Koszul as a commutative algebra if and only if it is Koszul as an associative algebra.
For $N-$Koszul algebras, there aren't any commutative examples (with more than one generator), since the commutativity relation is of degree two.

\begin{question}
Is there an operadic notion of $N-$Koszul algebra over a Koszul operad, with the property that a commutative algebra with all its generetors in degree $N$ is $N-$Koszul if and only if it satisfies the Cassidy-Shelton condition?
\end{question}

If there is such a condition, it might be formulated in terms of operadic cohomology, as the $N-$Koszul condition for associatice algebras can be formulated in terms of Hochschild cohomolgy.
The operadic cohomology for commutative algebras is the Harrison cohomology.

\begin{question}
What is the Harrison cohomology of coordinate ring of the generic determinantal variety?
\end{question}


\begin{thebibliography} {99}
\bibitem{AS87}M. Artin and W. Schelter, Graded algebras of global
  dimension 3, Adv. Math. 66 (1987), 171--216.
\bibitem{BF} Backelin, J.; Fr\"{o}berg, R.
Koszul algebras, Veronese subrings and rings with linear resolutions.
Rev. Roumaine Math. Pures Appl. 30 (1985), no. 2, 85--97. 
\bibitem{Bar} V. Baranovsky, BGG Correspondence for Projective Complete Intersections, Int. Math. Res. Not 45 (2005), 2759 --2774.
\bibitem{B01} R. Berger, 
Koszulity for nonquadratic algebras.
J. Algebra 239 (2001), no. 2, 705--734. 
\bibitem{BG} R. Berger and  V. Ginzburg,
Higher symplectic reflection algebras and non-homogeneous $N$-Koszul property.
J. Algebra 304 (2006), no. 1, 577--601. 
\bibitem{CS} T. Cassidy and B. Shelton, Generalizing the notion of Koszul algebra.
Math. Z. 260 (2008), no. 1, 93--114. 
\bibitem{Eis} D. Eisenbud, Commutative algebra.
With a view toward algebraic geometry. Graduate Texts in Mathematics, 150. Springer-Verlag, New York, 1995.
\bibitem{FV06} G. Fl\o ystad and J.E. Vatne, PBW-deformations of $N$-Koszul algebras. J. Algebra 302 (2006), no. 1, 116--155. 
\bibitem{GK} V. Ginzburg and M. Kapranov, Koszul duality for operads. Duke Math. J. 76 (1994), no. 1, 203--272.
\bibitem{GKR} A.L. Gorodentsev, A. Khoroshkin and A.N. Rudakov, On syzygies of highest weight orbits, Moscow seminar on Mathematical Physics. II, 79120, Amer.Math.Soc.Transl. Ser. 2, 221 (2007).
\bibitem{GMMVZ} E. L. Green, E.N. Marcos, R.  Mart\'{\i}nez-Villa and P. Zhang, $D$-Koszul algebras. J. Pure Appl. Algebra 193 (2004), no. 1-3, 141--162.
\bibitem{Lod1} J.-L. Loday, Cyclic homology. Grundlehren der matematischen Wissenschaften, 301 (2nd Edition), 1998. 
\bibitem{Lod2} J.-L. Loday, A. Frabetti, F. Chapoton and F. Goichot,Dialgebras and Related Operads, Springer Lecture Notes in mathematics, 1763, 2001.
\bibitem{M2} D.R. Grayson and M.E. Stillman, Macaulay2, a software system for research in algebraic geometry, Available at http://www.math.uiuc.edu/Macaulay2/
\bibitem{PP} A. Polishchuk and L. Positselski, Quadratic Algebras, AMS University Lecture series, 37, 2005.
\bibitem{Prid} S.B. Priddy, Koszul resolutions.
Trans. Amer. Math. Soc. 152 1970 39--60.
\bibitem{Val} B. Vallette, Manin products, Koszul duality, Loday algebras and Deligne conjecture
      Journal f\"{u}r die reine und angewandte Mathematik (Crelles Journal) Issue 620 (2008), pages 105--164
\end{thebibliography}
\end{document}